\documentclass{gtmon_a}
\pdfoutput=1

\usepackage{xy}        
\xyoption{all}
\usepackage{pinlabel}

\proceedingstitle{The Zieschang Gedenkschrift}
\conferencestart{5 September 2007}
\conferenceend{8 September 2007}
\conferencename{Conference in honour of Heiner Zieschang}
\conferencelocation{Toulouse, France}

\editor{Michel Boileau}
\givenname{Michel}
\surname{Boileau}

\editor{Martin Scharlemann}
\givenname{Martin}
\surname{Scharlemann}

\editor{Richard Weidmann}
\givenname{Richard}
\surname{Weidmann}

\title[Homology of the mapping class group $\Gamma_{2,1}$]{Homology of
the mapping class group $\Gamma_{2,1}$ for surfaces of genus 2 with a
boundary curve}

\author{Jochen Abhau}
\givenname{Jochen}
\surname{Abhau}

\author{Carl-Friedrich B\"odigheimer\newline Ralf Ehrenfried\newline }
\givenname{Carl-Friedrich}
\surname{B\"odigheimer}
\address{Mathematisches Institut\\
Universit\"at Bonn\\\newline
Beringstra{\ss}e 1\\
D--53115 Bonn\\
Germany}
\email{boedigheimer@math.uni-bonn.de}
\urladdr{}

\givenname{Ralf}
\surname{Ehrenfried}

\dedicatory{Dedicated to the memory of Heiner Zieschang}

\volumenumber{14}
\issuenumber{}
\publicationyear{2008}
\papernumber{1}
\startpage{1}
\endpage{25}

\doi{}
\MR{}
\Zbl{}

\arxivreference{0712.4254}


\received{15 October 2006}
\revised{21 March 2007}
\accepted{26 March 2007}
\published{29 April 2008}
\publishedonline{29 April 2008}
\proposed{}
\seconded{}
\corresponding{}
\version{}

\let\xysavmatrix\xymatrix
\def\xymatrix{\disablesubscriptcorrection\xysavmatrix}
\AtBeginDocument{\let\bar\wbar\let\tilde\wtilde\let\hat\what}

\def\strut{{\vrule height1.2em depth0.2em width0pt}}

\numberwithin{equation}{section}

\newcommand{\eins}{\mathbb{1}}

\newcommand{\del}{\partial}
\newcommand{\tdel}{\tilde{\partial}}
\newcommand{\homotopic}{\simeq}
\newcommand{\isomorphic}{\cong}

\newcommand{\norm}{\operatorname{norm}}

\newcommand{\ncyc}{\operatorname{ncyc}}
\newcommand{\sign}{\operatorname{sign}}

\newcommand{\perm}[1]{\langle #1 \rangle}

\newcommand{\Mod}{\mathfrak{Mod}}
\newcommand{\tMod}{\smash{\wwtilde{\mathfrak{M}}\mathfrak{od}}}

\newcommand{\Teich}{\mathfrak{Teich}}
\newcommand{\tTeich}{\smash{\wwtilde{\mathfrak{T}}\mathfrak{eich}}}

\newcommand{\Dip}{\mathfrak{Harm}}
\newcommand{\tDip}{\smash{\wwtilde{\mathfrak{H}}\mathfrak{arm}}}

\newcommand{\Par}{\mathrm{Par}}
\newcommand{\tPar}{\smash{\wwtilde{\mathrm{P}}\mathrm{ar}}}
\newcommand{\bPar}{\mathbb{P}}

\newcommand{\tGamma}{\wtilde{\Gamma}}

\newcommand{\tbPar}{\wtilde{\mathbb{P}}}

\newcommand{\tbQ}{\wwtilde{\mathbb{Q}}}

\newcommand{\bubu}{{\bullet\bullet}}
\newcommand{\topcell}{\mathrm{top}}

\newcommand{\mb}[1]{\mathbb{#1}}

\newcommand{\bC}{\mathbb{C}}

\newcommand{\bQ}{\mathbb{Q}}
\newcommand{\bR}{\mathbb{R}}
\newcommand{\bS}{\mathbb{S}}

\newcommand{\bZ}{\mathbb{Z}}

\newcommand{\cF}{\mathcal{F}}

\newcommand{\cK}{\mathcal{K}}

\newcommand{\cO}{\mathcal{O}}
\newcommand{\cP}{\mathcal{P}}
\newcommand{\cQ}{\mathcal{Q}}

\newcommand{\cT}{\mathcal{T}}

\newcommand{\cX}{\mathcal{X}}

\newcommand{\mf}[1]{\mathfrak{#1}}

\newcommand{\fS}{\mathfrak{S}}

\begin{document}

\begin{htmlabstract}
We report on the computation of the
integral homology of the mapping class group &Gamma;<sub>g,1</sub><sup>m</sup>
of genus g surfaces with one boundary curve and m punctures,
when 2g + m &le; 5, in particular &Gamma;<sub>2,1</sub><sup>0</sup>.
\end{htmlabstract}

\begin{abstract}
We report on the computation of the 
integral homology of the mapping class group $\Gamma_{g,1}^m$ 
of genus $g$ surfaces with one boundary curve and $m$ punctures, 
when $2g{+}m{\leq}5$, in particular $\Gamma_{2,1}^0$.
\end{abstract}

\begin{asciiabstract}
We report on the computation of the 
integral homology of the mapping class group Gamma_{g,1}^m
of genus g surfaces with one boundary curve and m punctures, 
when 2g + m < 6, in particular Gamma_{2,1}^0.
\end{asciiabstract}

\maketitle

\section{Introduction}
\label{sec1}

Let $\Mod = \Mod_{g,n}^m$ denote the moduli space of conformal equivalence 
classes of Riemann surfaces $F = F_{g,n}^m$ of genus $g \geq 0$ 
with $n \geq 1$ boundary curves and $m \geq 0$ permutable punctures.
One obtains an $(m!)$--fold covering space $\tMod = \tMod_{g,n}^m$ of $\Mod = \Mod_{g,n}^m$
if the punctures are declared not permutable.

Likewise, let $\Gamma = \Gamma_{g,n}^m$ be the corresponding mapping class group
of isotopy classes of orientation-preserving diffeomorphisms fixing the
boundary pointwise and possibly permuting the punctures.
The mapping class group where the punctures are to be fixed is a subgroup of 
index $m!$ in $\Gamma$ and is denoted by $\tGamma = \tGamma_{g,n}^m$.

The main result of this article is the computation of the integral homology
$H_*(\Gamma_{2,1}^0;\bZ)$  and $H_*(\Gamma_{2,1}^1;\bZ)$ 
of the mapping mapping class group for surfaces of genus $2$ 
with one boundary curve and with no (respectively, one) puncture.
These computations were done some years ago by the third author in \cite{Ehrenfried} 
and redone by the first author in \cite{Abhau},
both based on work by the second author in \cite{Boedigheimer1990} and \cite{Boedigheimer2005}. 
We give this belated report, because very few homology groups of mapping class groups are known;
prior to our computations the integral homology was known only for the easy case $g=1$.
See \fullref{sec2} for more remarks on previously known results.

$\Gamma$ and $\tGamma$ are torsion-free, since the diffeomorphisms are 
fixing at least one boundary curve.
They act therefore freely on the contractible Teichm\"uller spaces 
$\Teich = \Teich_{g,n}^m$ (respectively, $\tTeich = \tTeich_{g,n}^m$).
It follows, that the space $\Mod$ is a non-compact, connected (topological) manifold 
of dimension $d = 6g-6+3n+2m$;
the manifold $\tMod$ is always orientable,
but $\Mod$ is orientable only in the cases $m=0$ or $m=1$.
Furthermore, they have the homotopy type of the classifying space of the corresponding
mapping class group,  
\begin{equation}
\xymatrix{
\tMod\, \ar[d]\ar@{=}[r] & \,\tTeich/\tGamma \,\ar[d]\ar[r]^-{\homotopic} & \, B\tGamma \ar[d]  \\
\Mod \,       \ar@{=}[r] & \,\Teich/\Gamma   \,      \ar[r]^-{\homotopic} & \, B\Gamma  
}
\end{equation}
The key point of our method is to use a new description of 
the moduli space $\Mod$ (respectively, $\tMod$);
this is the space of parallel slit domains, 
the space of configurations of $h$ pairs of parallel, 
semi-infinite slits in $n$ complex planes;
here $h = 2g + m + 2n - 2$.

More precisely, there is a vector bundle 
$\Dip \to \Mod$ (respectively, $\tDip \to \tMod$) of dimension $d^* = m + 2n$,
and its one-point-compactification will be described by 
parallel slit domains.

To describe these bundles we first replace the moduli spaces above 
by the moduli space of closed Riemann surfaces $F$ of genus $g$, 
on which $n$ so-called dipole points $Q_1, \ldots, Q_n$ 
with non-zero tangent vectors $X_1, \ldots, X_n$, 
and $m$ permutable (respectively, non-permutable) punctures 
$P_1, \ldots, P_m$ are specified.
A conformal equivalence class is given by $[F,\cQ,\cP,\cX]$,
with $\cQ = (Q_1, \ldots, Q_n)$, $\cP = (P_1, \ldots, P_m)$,
and $\cX = (X_1, \ldots, X_n)$.
These new moduli spaces have the same homotopy type as $\Mod$
(respectively, $\tMod$)
or, in other words, the mapping class groups are isomorphic.
Note that their dimension is $d = 6g - 6 + 4n + 2m$,
since for the sake of simplicity we specified tangent vectors 
and not merely tangent directions. 

In the vector bundle $\Dip \to \Mod$ (respectively, $\tDip \to \tMod$)
the fibre over a point $[F,\cQ,\cP,\cX]$ in $\Mod$ (respectively, in
$\tMod$) consists of all harmonic functions 
$u \co F \to \bar{\bR} = \bR \cup \infty$ 
with a simple pole at each $Q_i$ having direction $X_i$ plus a logarithmic singularity, 
and with a logarithmic singularity at each $P_j$.
The bundle $\Dip$ is flat, more precisely 
$\Dip \isomorphic (\tMod \times_{\fS_{m}} \bR^{m}) \times \bR^{3n}$,
and the bundle $\tDip$ is trivial, being the pullback
of $\Dip$ along the covering $\tMod \to \tMod / \fS_{m} = \Mod$,
where $\fS_m$ is the $m$th symmetric group.

In \cite{Boedigheimer2005} we have
(following an earlier version in \cite{Boedigheimer1990})
introduced a finite cell complex $\Par = \Par(h,m,n)$
of dimension $d + d^* = 3h$
as a compactification of $\Dip$.
The complement of $\Dip$ is a subcomplex $\Par' \subset \Par$ of codimension $1$ and 
its points represent degenerate surfaces. 
A similar statement holds for $\tDip$.
The main result in \cite{Boedigheimer2005} and Ebert \cite{Ebert2005} is
\begin{eqnarray}
\Par - \Par'   \quad \isomorphic \quad & \Dip &  \quad \homotopic \quad \Mod ,  \\
\tPar - \tPar' \quad \isomorphic \quad & \tDip & \quad \homotopic \quad \tMod .
\end{eqnarray}
We later give some details when we exhibit the cellular chain complexes.

The pairs $(\Par,\Par')$ (respectively, $(\tPar,\tPar')$) are relative manifolds, 
and via Poincar{\'e} duality we obtain the isomorphisms: 
\begin{eqnarray}
H^*(\Par,\Par';\cO)     & \isomorphic &  H_{3h-*}(\Mod;\bZ) ,        \\
H^*(\tPar,\tPar';\bZ)   & \isomorphic &  H_{3h-*}(\tMod;\bZ)     
\end{eqnarray}            
These isomorphisms are the cap-product with a (relative) fundamental or orientation class
$[\mu]$ in $H^{3d}(\Par,\Par';\cO)$ (respectively, in $H^{3d}(\tPar,\tPar';\bZ)$)
where $\cO$ is the local coefficient system induced by the orientation covering.

In the following sections 
we concentrate for the sake of simplicity 
on the case of a single boundary curve ($n = 1$),
although the computations can be done in the general case.
To compute the homology of $\Mod$ and $\tMod$ 
we actually computed the cellular homology of the pairs $(\Par,\Par')$
with integral, with $\mod 2$ and with rational coefficients, 
and with coefficients in the orientation system $\cO$. 
The actual calculations were done using a computer program;
the program is not tracking the generators of the homology groups,
but some generators can easily be determined.

In \fullref{sec2} we will list our results.
We will give comments on previous and similar homology computations.
We give a description of the spaces $\Dip$ and $\tDip$ in \fullref{sec3},
including figures of configurations, that is, parallel slit domains;
in particular we give a complete list of cells for the 
case $g=1$ and $m=0$, which we will use in \fullref{sec6} for a demonstration.
In \fullref{sec4} we describe the cellular chain complexes.
Some properties of these chain complexes are described in \fullref{sec5}.
The calculation for the case $g=1$ and $m=0$ is done by hand in \fullref{sec6} 
to demonstrate the method, in particular the spectral sequence used. 
For the fundamental class $\mu$ and for some other homology classes 
we can give formulas (for their Poincar{\'e}) duals in \fullref{sec7}.
The orientation system is described in \fullref{sec8}.
And we comment on the computer program in \fullref{sec9}.

\subsection*{Acknowledgements} 
The authors are grateful to Jens Franke and Birgit Richter for valuable
discussions; thanks are also due to Marc Alexander Schweitzer and
Michael Griebel for access to computing capacities.  The referee made
several valuable suggestions and pointed out that some figures (see
\fullref{sec3}) will be helpful.

\section{Results}
\label{sec2}

\subsection{The permutable case}

We present our results in the form of a table.
It shows all non-trivial integral homology groups
of the moduli space $\Mod$ for $h \leq 5$ and $g > 0$ and $n = 1$.

\begin{center}
{\small
\begin{tabular}{|r|c|c|c|c|c|c|c|} \hline
\strut $* =$                         & 0              
                              & 1         
                              & 2              
                              & 3
                              & 4              
                              & 5
                              & 6               
                              \\ \hline
\strut $H_*(\Mod_{1,1}^0;\bZ)=$      & $\mb{Z}$       
                             & $\mb{Z}$ 
                             &                
                             &       
                             &
                             &
                             &
                             \\[1ex]
%
$H_*(\Mod_{1,1}^1;\bZ)=$       & $\mb{Z}$        
                              & $\mb{Z}$ 
                              & $\mb{Z}/2$     
                              &       
                              &                
                              & 
                              &                 
                              \\[1ex]
%
$H_*(\Mod_{1,1}^2,\bZ)=$       & $\mb{Z}$       
                              & $\mb{Z}{\oplus}\mb{Z}/2$ 
                              & $\mb{Z}/2{\oplus}\mb{Z}/2$  
                              & $\mb{Z}/2$      
                              &                 
                              &
                              &                
                              \\[1ex]
%
$H_*(\Mod_{1,1}^3;\bZ)=$       & $\mb{Z}$ 
                              & $\mb{Z}{\oplus}\mb{Z}/2$ 
                              & $\mb{Z}/2{\oplus}\mb{Z}/2$ 
                              & $\mb{Z}{\oplus}\mb{Z}/2{\oplus}\mb{Z}/2$       
                              & $\mb{Z}$
                              & $\mb{Z}$
                              &          
                              \\[1ex]
%
$H_*(\Mod_{2,1}^0;\bZ)=$       & $\mb{Z}$ 
                              & $\mb{Z}/10$ 
                              & $\mb{Z}/2$
                              & $\mb{Z}{\oplus}\mb{Z}/2$      
                              & $\mb{Z}/6$
                              &
                              &          
                              \\[1ex]
%
$H_*(\Mod_{2,1}^1;\bZ)=$       & $\mb{Z}$ 
                              & $\mb{Z}/10$ 
                              & $\mb{Z}{\oplus}\mb{Z}/2$
                              & $\mb{Z}{\oplus}\mb{Z}{\oplus}\mb{Z}/2{\oplus}\mb{Z}/2$      
                              & $\mb{Z}/6{\oplus}\mb{Z}/6$
                              & $\mb{Z}$
                              & $\mb{Z}$
                              \\ \hline
%
\end{tabular}
}
\end{center}

\subsection{Remarks}

(1)\qua
For $g = 0$, the moduli spaces $\Mod_{0,1}^m$ are 
the unordered configuration spaces of $m$ points in a disk or in the plane. 
The mapping class group $\Gamma_{0,1}^m$ is the (Artin) braid group 
$\mathrm{Br}(m)$ on $m$ strings.
Their homology is well-known and therefore we did not include the five 
cases $g=0$, $h=0, \ldots, 5$ in the table above.
See Arnol'd \cite{Arnold-1,Arnold-2}, Fuks \cite{Fuks}, 
Va\u{\i}n\v{s}te\u{\i}n \cite{Vainshtein} and Cohen \cite{Cohen-1,Cohen-2}. 

(2)\qua
The moduli space $\Mod_{1,1}^0$ of the torus with one boundary curve 
is homotopy equivalent to the complement of the trefoil knot in $\bS^3$.
The mapping class group $\Gamma_{1,1}^0$ is a central extension
\begin{equation}
1 \to \bZ \to \Gamma_{1,1}^0 \to \mathrm{SL}_2(\bZ) \to 1
\end{equation}
of $\mathrm{SL}_2(\bZ)$ by the integers 
and isomorphic to the third braid group $\mathrm{Br}(3)$.
This central extension is a special case of the central extension
\begin{equation}
1 \to \bZ^n \to \tGamma_{g,n}^m \to \tGamma_{g,0}^{n+m} \to 1 .
\end{equation}
Here the boundary curve is turned into a puncture by gluing-in a disc.
The connection to the sequence before 
is disguised by the fact $\tGamma_{1,0}^0 = \Gamma_{1,0}^0 = \Gamma_{1,0}^1$,
because a torus is an abelian variety; in this special case 
we have $\Gamma_{1,0}^0 = \Gamma_{1,0}^1 = \mathrm{SL}_2(\bZ)$.

(3)\qua
The computations of $H_*(\Mod_{2,1}^0)$ were done by the third author
\cite{Ehrenfried} some years ago, and redone by the first author in
\cite{Abhau}.  Very recently, Godin \cite{Godin} has obtained these
results independently, using a different method of computation.

(4)\qua
Note that very little is known about the (integral) homology 
of a single mapping class group. 
\begin{enumerate}
\item[(a)] Known are the first homology groups $H_1(\Gamma_{g,n}^m)$ for all values of $g,n,m$.
Apart from the classical case $g=1$,
this was proved by Mumford in \cite{Mumford} ($g=2, n=m=0$) 
and by Powell in \cite{Powell} ($g \geq 3, n=m=0$), 
which both use representations of the mapping class groups.
For the other values of $n$ and $m$ one can use the extensions 
\begin{eqnarray}
&1 \to \bZ^n \to \Gamma_{g,n}^m \to \Gamma_{g,0}^{n+m} \to 1  &
\label{eqn2.3}\\
&1 \to \mathrm{Br}_k(F_{g,n}^m) \to \Gamma_{g,n}^m \to \Gamma_{g,n}^{m+k}
\to 1  &\\
&1 \to \tGamma_{g,n}^m \to \Gamma_{g,n}^{m} \to \fS_m \to 1 &
\end{eqnarray}
or argue again via representations.
In the second extension $\mathrm{Br}_k(F) = \pi_1(C^k(F))$
is the $k$th braid group of the surface $F$, 
the fundamental group of the $k$th unordered configuration space $C^k(F)$.
See Korkmaz \cite{Korkmaz-2} for a survey.

\item[(b)] Also known are the second homology groups $H_2(\Gamma_{g,n}^m)$ 
by work of Harer \cite{Harer-H2};
see also the corrections by Harer \cite{Harer-stability} and Morita \cite{Morita}.
These computations are based on the curve complex.
For a proof using presentations see Pitsch \cite{Pitsch}
and Korkmaz--Stipsicz \cite{Korkmaz-Stipsicz}. 

\item[(c)] For the third homology group $H_3(\Gamma_{g,n}^m; \mb{Q})$ 
the reference is Harer \cite{Harer-H3}.
\end{enumerate}

(5)\qua
The stable mapping class groups $\Gamma_{\infty,n}^m$,
that is, the limit of the mapping class groups 
for the genus tending to infinity,
can be defined for surfaces with boundary curves:
gluing on a torus with two boundary curves embeds 
$\Gamma_{g,n}^m$ into $\Gamma_{g+1,n}^m$. 
This embedding induces a homology isomorphism
in degrees $\leq (g-1)/3$. 
This is Harer's stability theorem,
see \cite{Harer-stability}, and \cite{Ivanov}.
Likewise,
gluing a disc onto a boundary curve and recording its center as a new puncture
gives a surjection of mapping class groups $\Gamma_{g,n+1}^m \to \Gamma_{g,n}$; 
see the extension \eqref{eqn2.3} above. 
This surjection induces also a homology isomorphism in
degrees $\leq g/3$; see \cite{Harer-survey}.

The proof of the Mumford conjecture by Madsen and Weiss
\cite{Madsen-Weiss} determines the homology of the stable mapping class
group with rational coefficients.  Using this, Galatius \cite{Galatius}
has computed the mod--$p$ homology of the stable mapping class group.

Note that none of our computations belong to the stable range.

(6)\qua
The mapping class groups for surfaces with (respectively, without) boundary
are quite different.  The former are torsion-free groups, and have
finite-dimensional classifying spaces, which can be chosen to be
manifolds.  Thus there homology is of finite type.  In contrast,
the mapping class groups of surfaces without boundary contain
elements of finite order.  Thus their classifying spaces can not
be finite-dimensional.  Nevertheless, they have finite virtual
homological dimension by a theorem of Harer \cite{Harer-homdim}.
In the homology of mapping class groups without boundary Glover--Mislin
\cite{Glover-Mislin} have found torsion elements by considering finite
subgroups.

(7)\qua
For homology computations with field coefficients see for example the
work of Benson and Cohen \cite{Benson-Cohen}, and the second author,
Cohen and Peim \cite{Boedigheimer-Cohen-Peim}, for $\Gamma_{2,0}^0$ with
coefficients in the field $\mathbb{F}_p$ and $p = 2,3,5$; or see the work
of Looijenga \cite{Looijenga} for $\Gamma_{3,0}^0$ and $\Gamma_{3,0}^1$,
and the work of Tommasi \cite{Tommasi-1,Tommasi-2} for $\Gamma_{3,0}^2$
and $\Gamma_{4,0}^0$, all with rational coefficients.

(8)\qua
The Hilbert uniformization method can also be used to parametrize 
the moduli spaces of surfaces with incoming and outgoing boundary curves, 
see \cite{Boedigheimer2003}.

(9)\qua
This uniformization method can also be used 
for moduli spaces of non-orientable surfaces (Kleinian surfaces); 
see the work of Ebert \cite{Ebert2003} and Zaw \cite{Zaw}.
One needs besides the {\it oriented} slit pairs 
(or handles) we use here 
a second kind, namely the {\it unoriented} slit pairs (or cross-caps). 
In \cite{Zaw} the resulting cell decomposition of the moduli space 
was used to calculate the homology of non-orientable 
mapping class groups; see also Korkmaz \cite{Korkmaz-1}.

\subsection{The non-permutable case}

Obviously, $\tMod_{g,n}^m = \Mod_{g,n}^m$ for $m = 0$ and $m = 1$.

For $g = 0$ the moduli space $\tMod_{0,1}^m$ 
is the $m$th ordered configuration space of a disk or a plane. 
And the mapping class group $\tGamma_{0,1}^m$ is the pure (Artin) braid group 
$\mathrm{Br}(m)$ on $m$ strings. 
As before, their homology is well-known (see for example Arnol'd \cite{Arnold-1}) 
and therefore we did not include the cases
$g = 0$ and $m = 0, \ldots, 5$ in this table.

Thus there are for $h \leq 5$ only the cases $g=1$, $m=2$ and
$g=1$, $m=3$.  With the available capacities we could only finish
the computations for the first case, which is listed below.
See \fullref{sec9} for comments.

\begin{center}
\begin{tabular}{|r|c|c|c|c|c|} \hline
\strut $* =$                         & 0              
                               & 1         
                                & 2              
                                 & 3
                                  & 4              
                                   \\ \hline
\strut $H_*(\tMod_{1,1}^2;\bZ) = $  &  $\mb{Z}$     
                              &  $\mb{Z}$  
                               &  $\mb{Z}/2\oplus\mb{Z}/2\oplus\mb{Z}/2$  
                                &  $\mb{Z}$      
                                 &  $\mb{Z}$               
                                  \\ \hline
\end{tabular}
\end{center}


\section{Parallel slit domains}
\label{sec3}

\subsection{Dissecting a surface}

We shall give here a brief description of the Hilbert uniformization,
thereby concentrating on the case $n=1$ and on the permutable case.

The Hilbert uniformization associates to a point in $\Dip$
a so-called parallel slit domain.
Consider on the surface $\cF = [F, \cQ, \cX, \cP]$ a harmonic function 
$u \co F \to \bar{\bR} = \bR \cup \infty$ 
with exactly one dipole singularity at the point $Q_1$
and logarithmic singularities at the punctures $P_1, \ldots, P_m$. 
More precisely, 
for a local chart $z$ with $z(Q_1) = 0$ and $dz(X_1) = 1$ 
we require that for some positive real number $A_1$ and $B_1$ the function   
$u(z) - A_1 \mathrm{Re}(1/z) - B_1 \mathrm{Re}(\log(z)) $ is harmonic,
and further that for a local chart $z$ with $z(P_i) = 0$ 
we require that the function 
$u(z) - C_i \mathrm{Re}(\log(z))$ is harmonic.
Obviously, the residues of these singularities must 
satisfy the equation $-B_1 + C_1 + \dots + C_m = 0$. 
The set of these auxiliary parameters is the fibre of the vector bundle 
$\Dip \to \Mod$ over the point $\cF$.
Note that the function $u$ is -- up to a constant -- uniquely determined 
by $\cF$ and the auxiliary parameters.

This function $u$ will have up to $h$ critical points $S_1, S_2, \ldots $,
where $u(z) = \mathrm{Re}(z^{l_k})$, for $l_k \geq 2$, in some local chart.
The sum of the indices $\mathrm{ind}(S_k)$ plus the index $-2$ at the dipole $Q_1$ 
must be the Euler characteristic $\chi = 2 - 2g + m$.

The stable critical graph $\cK$ consists of the vertices $Q_1$ and $S_1, S_2, \ldots$, 
and all gradient curves of $u$ leaving critical points 
(to either go to $Q_1$ or to another critical point) as edges.
Its complement $F_0 = F - \cK$ is simply-connected.
Thus $u$ is on $F_0$ the real part of 
a holomorphic map $w = u + \sqrt{-1} v \co F_0 \to \bC$. 
Note that $w$ is unique up to a constant.

\subsection{Parallel slit domains}

The image of $w$ is the complex plane cut along horizontal slits
starting at some point and running towards infinity on the left.
This is called a {\it parallel slit domain} in geometric function theory.
These slits come - generically - in pairs, suitably interlocked.

To see this, let us assume that $u$ is generic, that is 
(1) that $u$ is locally a Morse function ($\mathrm{ind}(S_1) = \mathrm{ind}(S_2) = \cdots = 1$),  
(2) that the critical values are all distinct, and 
(3) that a gradient curve leaving a critical point does not enter another critical point,
but runs to the dipole $Q_1$.

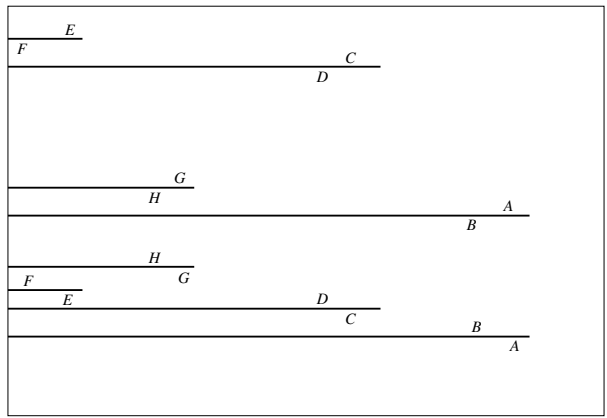
\begin{figure}[ht!]
\begin{center}
\labellist\tiny
\hair1.5pt
\pinlabel {$E$} [b] at 62 366
\pinlabel {$F$} [t] at 16 365
\pinlabel {$C$} [b] at 333 339
\pinlabel {$D$} [t] at 306 338
\pinlabel {$G$} [b] at 168 222
\pinlabel {$H$} [t] at 144 221
\pinlabel {$A$} [b] at 486 195
\pinlabel {$B$} [t] at 450 194
\pinlabel {$H$} [b] at 144 145
\pinlabel {$G$} [t] at 172 144
\pinlabel {$F$} [b] at 22 123
\pinlabel {$E$} [t] at 60 122
\pinlabel {$D$} [b] at 306 105
\pinlabel {$C$} [t] at 333 104
\pinlabel {$B$} [b] at 455 78
\pinlabel {$A$} [t] at 492 77
\endlabellist
\includegraphics[scale=0.4]{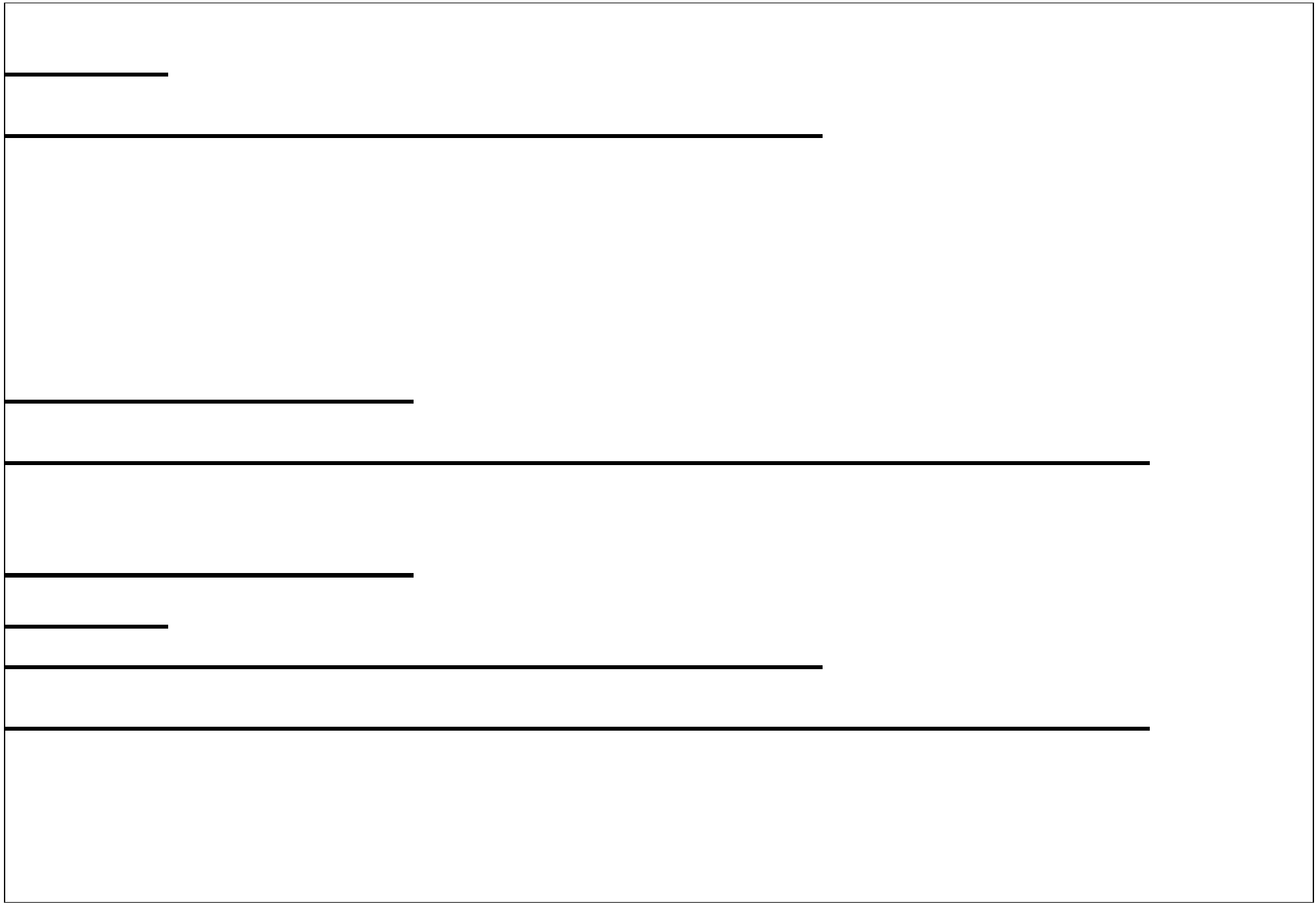}
\caption{An example of a generic parallel slit domain 
with $g=2$, $m=0$ and $n=1$}
\label{fig1}
\end{center}
\end{figure}

For $g=1, n=1, m=0$ there are two generic configurations;
see the following figures.
For $g=2, n=1, m=0$ there are $504$ generic configurations;
\fullref{fig1} shows such an example.
The letters $A, \ldots, H$ indicate the boundary identifications: 
the right bank of a slit is glued to the left bank of the paired slit.
The function $u$ corresponds to the real part.
The dipole point corresponds to infinity.

\begin{figure}[ht!]
\begin{center}
\labellist\tiny
\hair2pt
\pinlabel {$C$} [b] at 200 214
\pinlabel {$D$} [t] at 168 214
\pinlabel {$A$} [b] at 95 123
\pinlabel {$B$} [t] at 67 123
\pinlabel {$D$} [b] at 163 96
\pinlabel {$C$} [t] at 204 96
\pinlabel {$B$} [b] at 68 56
\pinlabel {$A$} [t] at 100 56
\pinlabel {$C$} [b] at 622 213
\pinlabel {$D$} [t] at 590 213
\pinlabel {$D$} [b] at 590 146
\pinlabel {$C$} [t] at 622 146
\pinlabel {$A$} [b] at 522 108
\pinlabel {$B$} [t] at 468 108
\pinlabel {$B$} [b] at 468 56
\pinlabel {$A$} [t] at 522 56
\endlabellist
\includegraphics[scale=0.5]{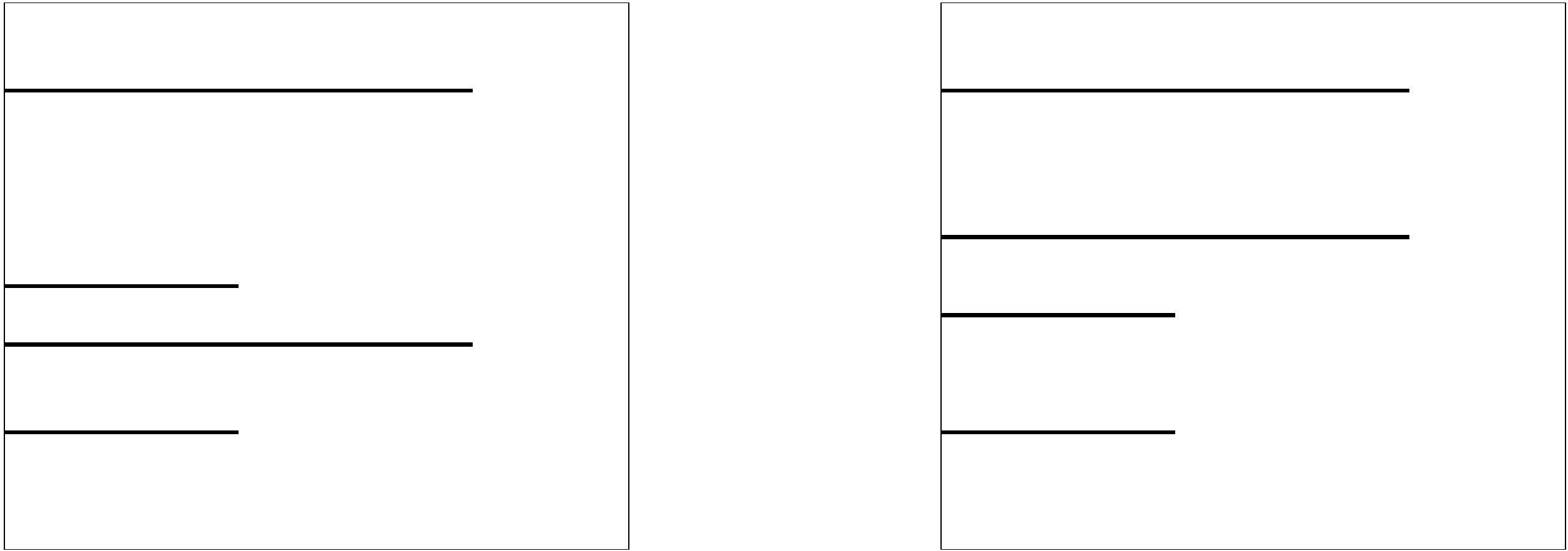}
\caption{Two examples of generic parallel slit domains 
with two pairs of slits: the different interlocking patterns 
give $g=1, m=0$ for the configuration on the left,
and $g=0, m=2$ for the configuration on the right}
\label{fig2}
\end{center}
\end{figure}
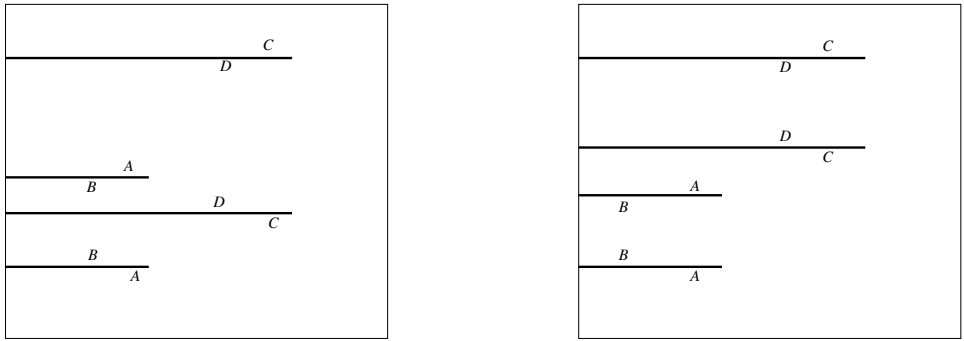

\fullref{fig2} shows two generic configurations with two pairs
of slits; but since the interlocking pattern is different, for the left
configuration $g=1$, $m=0$, whereas for the right one $g=0$, $m=2$.

\subsection{Configurations and cells}
\label{sec3.3}

If we identify the real line with an interval,
the vertical and horizontal distances between the slits 
are two systems of (barycentric) coordinates.
In other words, a top-dimensional cell is a 
product $\Delta^{2h} \times \Delta^h$ of two simplices,
where $h = 2g+m$ is the number of pairs of slits.

In the non-generic case, the slits may have equal length,
and may touch each other; 
in the latter case a shorter slit can jump to the other side of the 
partner of the longer slit pair. 
See the rows in \fullref{fig6} for examples of such jumps.
This leads to identifications among faces of top-dimensional cells.
Furthermore, not all configurations are allowed,
since they may lead to surfaces with singularities or to surfaces with too
small a genus or too few punctures.
Both types of cells are called degenerate and will be excluded from the space $\Dip$.

We refer the reader to \cite{Boedigheimer2005} for details.

\begin{figure}[ht!]
\begin{center}
\labellist\small
\hair2pt
\pinlabel {$B$} [t] at 19 186
\pinlabel {$A$} [b] at 90 186
\pinlabel {$A$} [t] at 90 96
\pinlabel {$B$} [b] at 17 96
\pinlabel {$20$} at 54 64
\pinlabel {$21$} at 54 113
\pinlabel {$22$} at 54 163
\pinlabel {$23$} at 54 208
\pinlabel {$24$} at 54 271
\pinlabel {$\sigma_2$} at 54 10
\pinlabel {$D$} [t] at 132 254
\pinlabel {$D$} [b] at 132 150
\pinlabel {$C$} [b] at 199 254
\pinlabel {$C$} [t] at 199 150
\pinlabel {$10$} at 163 64
\pinlabel {$11$} at 163 109
\pinlabel {$12$} at 163 163
\pinlabel {$13$} at 163 203
\pinlabel {$14$} at 163 275
\pinlabel {$\sigma_1$} at 163 10
\pinlabel {$00$} at 249 64
\pinlabel {$01$} at 249 109
\pinlabel {$02$} at 249 163
\pinlabel {$03$} at 249 208
\pinlabel {$04$} at 249 271
\pinlabel {$\sigma_0$} at 249 10
\pinlabel {$B$} [t] at 450 150
\pinlabel {$B$} [b] at 450 97
\pinlabel {$A$} [t] at 522 97
\pinlabel {$A$} [b] at 522 150
\pinlabel {$20$} at 487 59
\pinlabel {$21$} at 487 113
\pinlabel {$22$} at 487 163
\pinlabel {$23$} at 487 212
\pinlabel {$24$} at 487 266
\pinlabel {$\sigma_2$} at 487 10
\pinlabel {$D$} [b] at 564 185
\pinlabel {$D$} [t] at 564 252
\pinlabel {$C$} [b] at 631 252
\pinlabel {$C$} [t] at 631 185
\pinlabel {$14$} at 591 266
\pinlabel {$10$} at 595 59
\pinlabel {$11$} at 595 113
\pinlabel {$12$} at 595 163
\pinlabel {$13$} at 595 212
\pinlabel {$\sigma_1$} at 595 10
\pinlabel {$00$} at 685 59
\pinlabel {$01$} at 685 113
\pinlabel {$02$} at 685 163
\pinlabel {$03$} at 685 212
\pinlabel {$04$} at 685 266
\pinlabel {$\sigma_0$} at 685 10
\endlabellist
\includegraphics[scale=0.5]{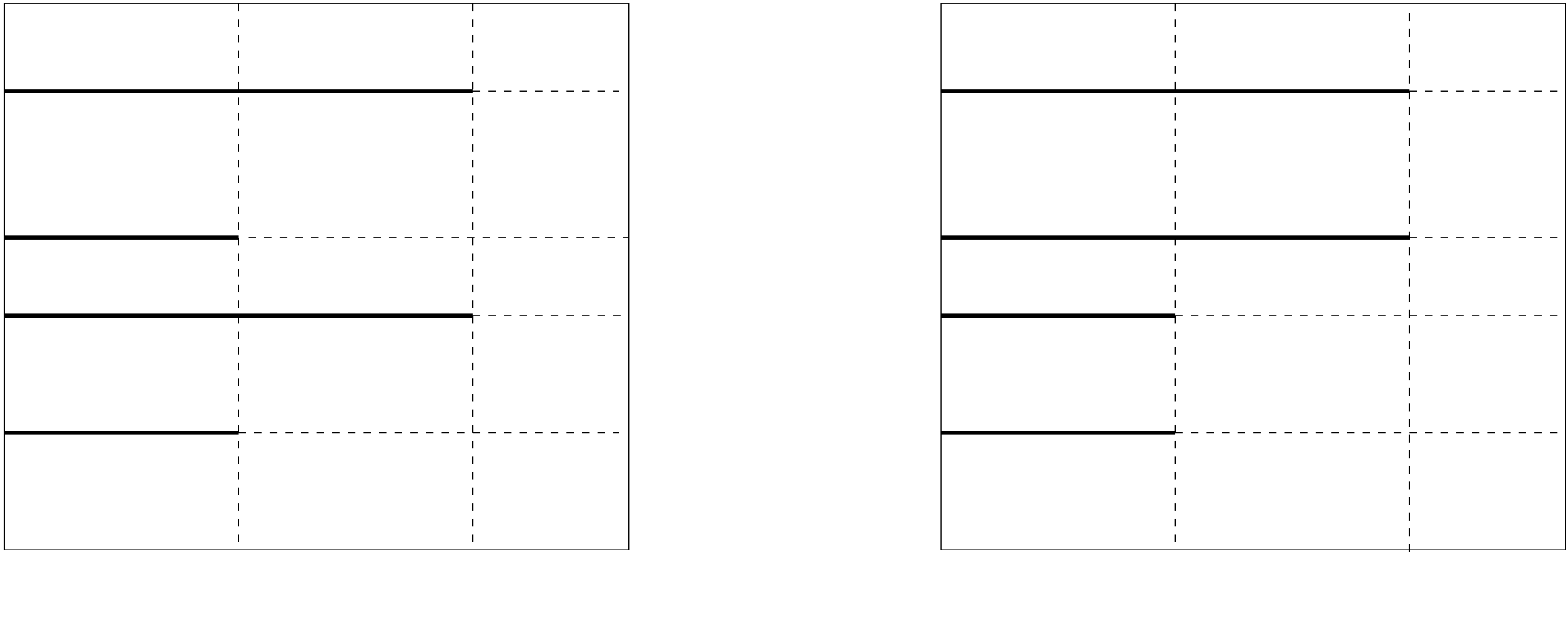}
\caption{The same two examples as in \fullref{fig2},
now shown with the grid pattern and the numbering of the rectangles;
the example on the left is the cell $\Sigma_1$}
\label{fig3}
\end{center}
\end{figure}

\fullref{fig3} shows how the combinatorial type of these two configurations
are codified. 
The horizontals of the slits and 
the verticals at the slit ends give a grid;
we number the columns $0,1, \ldots, q\leq2$ from right to left,
and the rows $0, 1, \ldots, p\leq4$ from bottom to top.
Note that in general $q \leq h$ and $p \leq 2h$.
Let $R_{i,j}$ denote the $j$th rectangle in the $i$th column.

For the $i$th column let $\sigma_i$ denote the permutation 
describing the re-gluing of the slit plane:
the upper edge of the rectangle $R_{i,j}$ 
is glued to the lower edge of $R_{i,\sigma_i(j)}$
(for $i = 0, 1, \ldots, q; \,\, j = 0, 1, \ldots, p-1$).
And the left edge of $R_{i,j}$ is glued to the right edge of $R_{i+1,j}$
(for $i = 0, 1, \ldots, q-1; \,\, j = 0, 1, \ldots, p$).

We write a permutation $\sigma$ as a product of its cycles, 
and a cycle sending $i_0$ to $i_1$, and $i_1$ to $i_2$, and so on, and $i_l$ to $i_0$,
is written from right to left, namely
$\perm{i_l, \ldots, i_2, i_1, i_0}$.  
For better readability we drop the commas between the one-digit numbers 
and write 
$\perm{i_l \, \ldots \, i_2 \, i_1, \, i_0}$.

In our examples in \fullref{fig3} these permutations are
\begin{align*}
\sigma_0 &= \perm{4\,3\,2\,1\,0}, &
\sigma_1 &= \perm{4\,1\,0}\perm{3\,2}, &
\sigma_2 &= \perm{4\,1\,2\,3\,0} \\
\text{and}\qquad \sigma_0 &= \perm{4\,3\,2\,1\,0}, &
\sigma_1 &= \perm{4\,2\,1\,0}\perm{3}, &
\sigma_2 &= \perm{4\,2\,0}\perm{1}\perm{3},
\end{align*}
respectively.

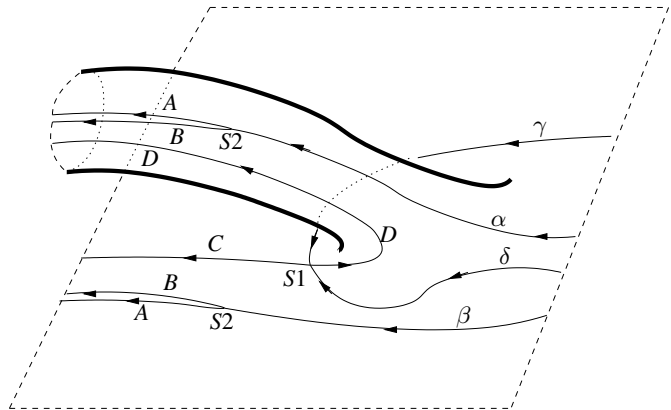
\begin{figure}[ht!]
\begin{center}
\labellist\footnotesize
\hair2pt
\pinlabel {$D$} [t] at 89 166
\pinlabel {$A$} [b] at 102 185
\pinlabel {$B$} [t] at 106 180
\pinlabel {$S2$} [t] at 140 177
\pinlabel {$A$} [t] at 84 69
\pinlabel {$B$} [b] at 102 72
\pinlabel {$S2$} [t] at 134 64
\pinlabel {$C$} [b] at 130 96
\pinlabel {$S1$} [tr] at 190 90
\pinlabel {$D$} [l] at 230 110
\pinlabel {$\beta$} [b] at 286 48
\pinlabel {$\alpha$} [b] at 308 112
\pinlabel {$\delta$} [b] at 312 88
\pinlabel {$\gamma$} [b] at 334 168
\endlabellist
\includegraphics[scale=0.6]{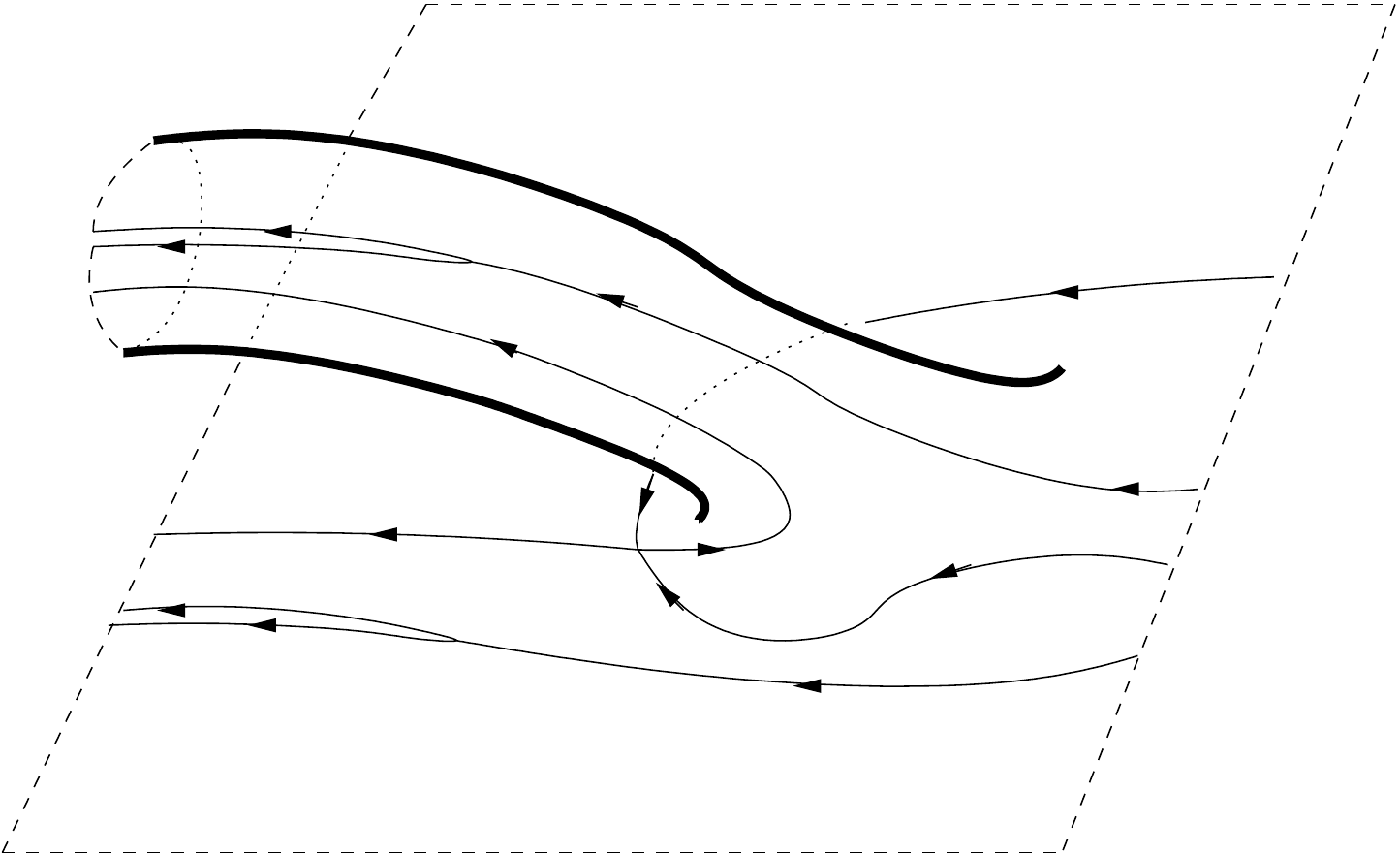}
\caption{This shows the gluing in process: the surface is half-way finished}
\label{fig4}
\end{center}
\end{figure}

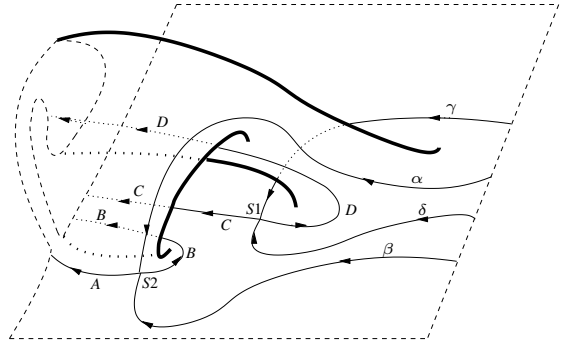
\begin{figure}[ht!]
\begin{center}
\labellist\tiny
\hair2pt
\pinlabel {$A$} [t] at 66 50
\pinlabel {$B$} [b] at 70 86
\pinlabel {$C$} [b] at 98 104
\pinlabel {$D$} [b] at 116 156
\pinlabel {$S2$} [tl] at 98 50
\pinlabel {$B$} [l] at 130 65
\pinlabel {$C$} [t] at 164 95
\pinlabel {$S1$} [br] at 194 92
\pinlabel {$D$} [l] at 250 98
\pinlabel {$\beta$} [b] at 286 58
\pinlabel {$\alpha$} [b] at 308 114
\pinlabel {$\delta$} [b] at 312 90
\pinlabel {$\gamma$} [b] at 334 166
\endlabellist
\includegraphics[scale=0.5]{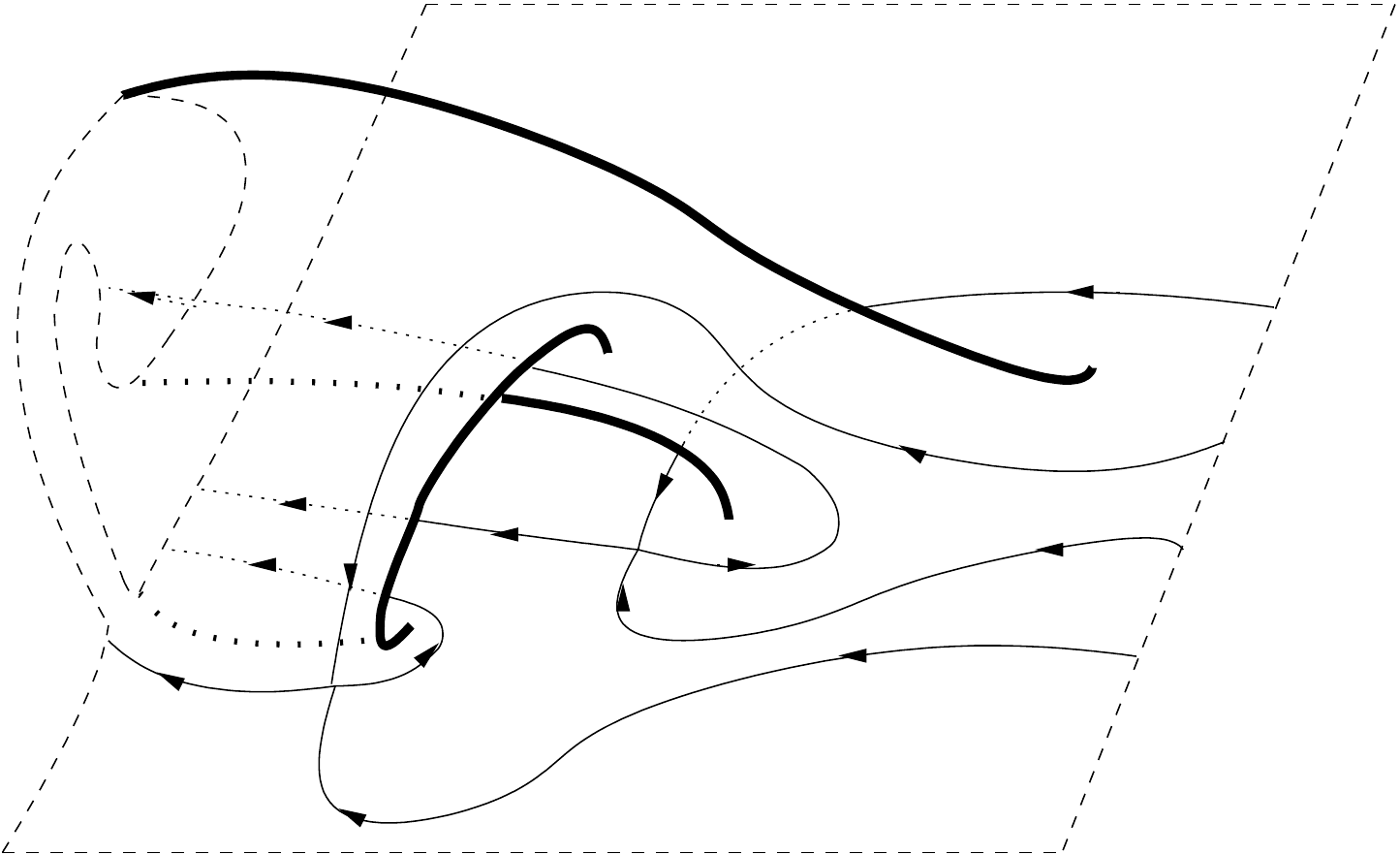}
\caption{The finished surface}
\label{fig5}
\end{center}
\end{figure}

\begin{figure}[ht!]
\begin{center}
\includegraphics[scale=0.6]{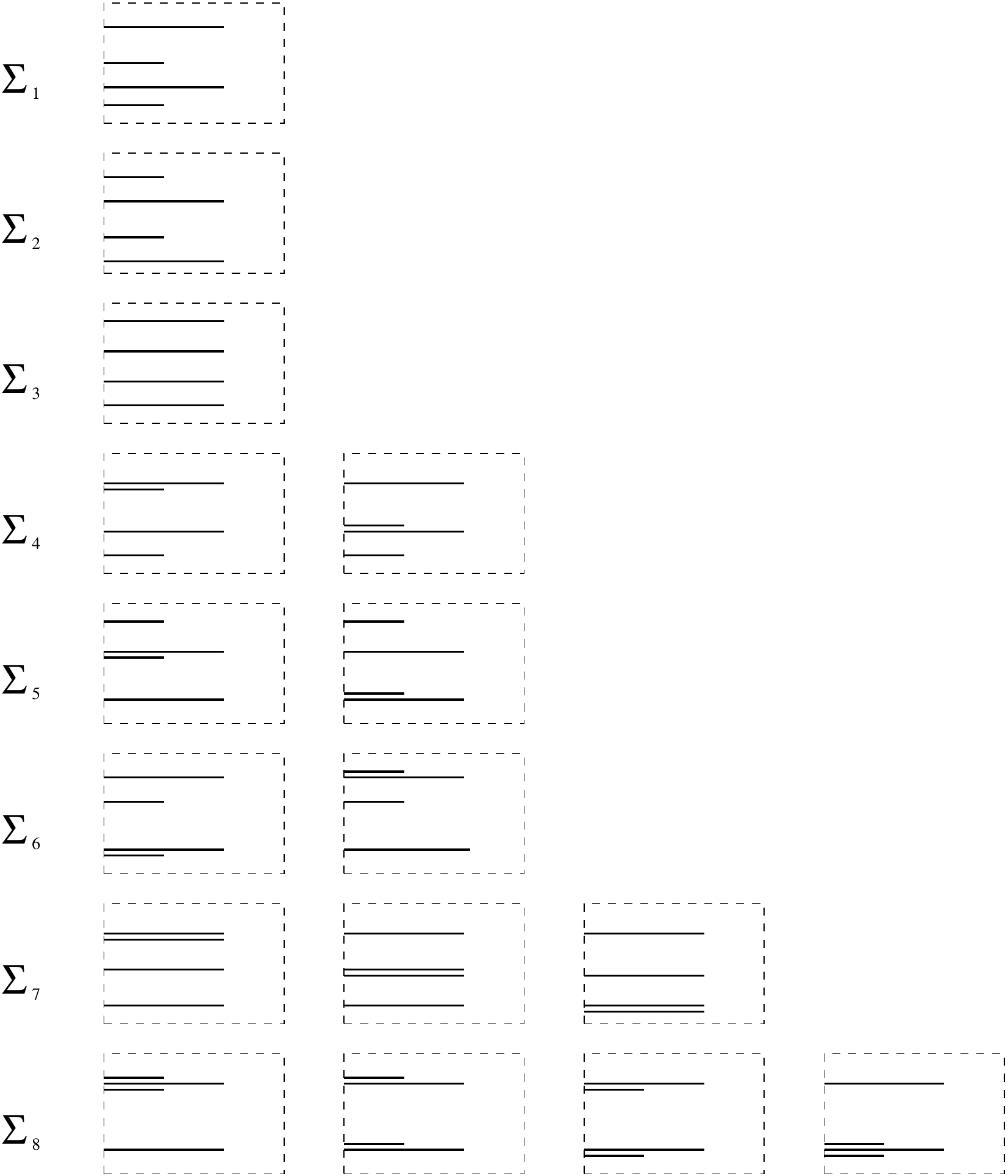}
\caption{All eight cells for $g=1$, $m=0$ and $n=1$:
in each row the configurations representing this cell are shown}
\label{fig6}
\end{center}
\end{figure}

Figures \ref{fig4} and \ref{fig5} show the gluing process
in two steps.  The example is the left hand side of \fullref{fig3}, a
generic configuration with $g=1$, $n=1$ and $m=0$, or $\Sigma_1$ from the
list in \fullref{sec3.4}.  In \fullref{fig4} the first slit pair is glued, producing
a tube growing out of the complex plane; the second slit pair is still
seen as two slits, one on the tube, one outside the tube.  The gradient
field of $u$ is shown, with with two critical points $S_1$ and $S_2$.

\fullref{fig5} shows the finished surface.
The tube is with its later part re-glued into the plane,
creating a plane with a handle.
The point at infinity needs to be added.

\subsection{The example $g=1$, $n=1$ and  $m=0$ geometrically}
\label{sec3.4}

Here we want to give a complete list of all non-degenerate cells
for the case $g=1$, $m=0$ and $n=1$; this illustrates the explanations above
and will be used later.

\fullref{fig4} shows all eight non-degenerate cells $\Sigma_1, \ldots, \Sigma_8$.
Recall that a cell is an equivalence class of configurations;
thus in each row we have drawn all configurations representing a cell.

Next we give the permutations of these eight cells. 
They will be used in \fullref{sec6} when we study this example algebraically.

For $(p,q) = (4,2)$ there are two non-degenerate cells
$$\Sigma_1 = (\perm{41230}, \perm{32}\perm{410}, \perm{43210} )
\quad\text{and}\quad
\Sigma_2 = (\perm{41230}, \perm{21}\perm{430}, \perm{43210} ).$$
For $(p,q) = (4,1)$ there is one non-degenerate cell
$$\Sigma_3 = (\perm{41230}, \perm{43210} ).$$
For $(p,q) = (3,2)$ there are three non-degenerate cells
\begin{align*}
\Sigma_4 &= (\perm{3120}, \perm{2}\perm{310}, \perm{3210} ), \\
\Sigma_5 &= (\perm{3120}, \perm{320}\perm{1}, \perm{3210} ), \\
\text{and}\quad \Sigma_6 &= (\perm{3120}, \perm{21}\perm{30}, \perm{3210} ).
\end{align*}
For $(p,q) = (3,1)$ there is one non-degenerate cell 
$$\Sigma_7 = (\perm{3120}, \perm{3210} ).$$
For $(p,q) = (2,2)$ there is one non-degenerate cell 
$$\Sigma_8 = (\perm{210}, \perm{20}\perm{1}, \perm{210} ).$$
Note that the number of cells grows fast with $g$. 
In the next case $g=2, n=1, m=0$ 
we have $17136$ non-degenerate cells in dimensions $5$ to $12$.

\section{The chain complexes  $\bPar_{\bubu}$ and $\tbPar_{\bubu}$}
\label{sec4}

Let $h$ and $m$ be fixed and set $n=1$.

\subsection{The cells}

The cellular chain complex of the space $\Par = \tPar(h,m,n)$ 
is the total complex of a bi-graded complex 
$\bPar_\bubu = \bigoplus \bPar_{p,q}$,
where $0 \leq p \leq 2h$ and $0 \leq q \leq h$.
Here $\bPar_{p,q}$ is the free abelian group generated by all
$(q{+}1)$--tuples $\Sigma = (\sigma_q, \ldots, \sigma_0)$ of permutations 
$\sigma_i$ in the symmetric group $\mf{S}_{p+1}$, 
acting on $[p] = \{0,1, \ldots, p \}$,
and satisfying the following conditions:
\begin{align}
\norm(\Sigma)  &\leq h \label{eqn4.1}\\
\ncyc(\sigma_q) &\leq  m+1  \label{eqn4.2} 
\end{align}
Here the $\norm(\Sigma)$ is the sum of the word lengths, 
$$\norm(\Sigma) = \ell(\sigma_q \sigma_{q-1}^{-1}) + \cdots +
\ell(\sigma_1 \, \sigma_0^{-1}),$$
where the word length $\ell(\alpha)$ is 
measured with respect to the generating set of all transpositions in $\mf{S}_{p+1}$.
And $\ncyc(\alpha)$ is the number of cycles of a permutation $\alpha$.

The cellular chain complex of the space $\tPar$ 
is the total complex of the bi-graded complex 
$\tbPar_\bubu$,
where $\tbPar_{p,q}$ is generated by elements 
$(\nu; \Sigma) = (\nu; (\sigma_q, \ldots, \sigma_0))$
where $\Sigma$ is as before a $(q{+}1)$--tuple of permutations
in $\fS_{p+1}$ satisfying the condition \eqref{eqn4.1} and \eqref{eqn4.2} above;
and $\nu \co [p] \to [m]$ is a function decoding a numbering of the punctures,
which therefore must satisfy the following conditions:
\begin{align}
& \nu \text{ is invariant under } \sigma_q \label{eqn4.3} \\
& \nu \text{ induces a bijection } [p]/{\sigma_q} \to [m] \label{eqn4.4} \\ 
& \nu(0) = \nu(p) = 0 \label{eqn4.5}
\end{align} 
Here $[p]/\sigma_q$ denotes the set of $\sigma_q$--orbits in $[p]$.

\subsection{The surface associated to a configuration}
\label{sec4.2}

We repeat here the general construction we described for a generic 
configuration in \fullref{sec3.3}

Given a cell $\Sigma$ of bi-degree $(p,q)$ 
and a point $a = (a_0, \ldots, a_p)$ in a $p$--dimensional simplex $\Delta^p$
and a point $b = (b_0, \ldots, b_q)$ in a $q$--dimensional simplex $\Delta^q$ 
we associate a surface $F = F(\Sigma; a,b)$ as follows.
To begin, we identify the open unit interval with the real line,
and thus interpret the $a_i$ as real and the $b_j$ as imaginary numbers.

Subdivide the complex plane by $p$ horizontal lines at the $a_i$
(respectively, by $q$ vertical lines at the $b_j$)
and number the $(p+1)(q+1)$ rectangles $R_{i,j}$ from bottom to top
(respectively, from right to left).  See the examples in \fullref{sec3.3}.

If the rectangles $R_{i,j}$ along a column $i$ are glued in cycles
according to the permutation $\sigma_i$, and along a row according to
their numbering, one obtains an ``open'' $2$--complex which we compactify
by adding a point $Q_1$ ``at infinity'' and at most $m$ points $P_1,
\ldots, P_m$ at ends of the tubes created by the cycles of rectangles
$R_{q,j}$; we denote this $2$--complex by $F = F(\Sigma; a,b)$.  Note that
the $(p-1)(q-1)$ rectangles $R_{i,j}$ with $0 < i < q$ and $0 < j < p$ are
now quadrilaterals on $F$, those with $i=0, 0 < j < p$, $i=q, 0 < j < p$,
$0 < i < q, j=0$ or $0 < i < q, j=p$ are triangles (with one corner being
$Q_1$) and the remaining four are two-gones with one corner being $Q_1$.

The definition of a degenerate configuration in \fullref{sec4.3} is made
such $F$ is a surface of genus exactly $g$ with exactly $m$ punctures
if and only if $\Sigma$ is non-degenerate.  (Note that we use the term
``degenerate'' also for smooth $F$ with genus or number of punctures
smaller than expected.)  Then we take $F$ with the obvious complex
structure and at the point $Q_1$ we take the direction $X_1$ corresponding
to the horizontal direction.  Note that there is a harmonic function $u$
given by projection to the real axis.

If $F$ is a surface, but the genus is less than $g$ or there are less
than $m$ punctures, then $\Sigma$ will also be called degenerate.

\subsection{Degenerate cells}
\label{sec4.3}

The subcomplex $\bPar'$, which corresponds to the subspace $\Par'$ of
$\Par$, consists of all $\Sigma$ with $\norm$ less than $h$, or with a
$\sigma_q$ having less than $m+1$ cycles, or with any of the following
conditions violated:
\begin{align}
& \sigma_i(p) = 0 \text{ for } i = 0, \ldots, q \label{eqn4.6}\\
& \sigma_0 \text{ is the rotation } \omega_p = \perm{p,\ldots, 1, 0}
\label{eqn4.7} \\
& \sigma_{i+1} \not= \sigma_i \text{ for } i= 0, \ldots, q \label{eqn4.8} \\
& \text{There is no }  0 \leq k \leq p-1  
  \text{ such that } \sigma_i(k) = k+1 \text{ for all } i=0, \ldots, q
\label{eqn4.9}
\end{align}
The subcomplex $\tbPar'$,
which corresponds to the subspace $\tPar'$,  
consists of all cells $(\nu;\Sigma)$ 
with $\norm$ less than $h$,
or with $\sigma_q$ having less than $m+1$ cycles, 
or with $\nu$ not surjective, 
or with any of the conditions \eqref{eqn4.3}--\eqref{eqn4.5} or
\eqref{eqn4.6}--\eqref{eqn4.9} violated.

Note that the first and the last faces of any cell under $\del'$ and $\del''$
are degenerate.

The cells in $\tPar'$ (respectively, in $\Par'$)
are called degenerate
although the surfaces represented by their points 
might be smooth surfaces -- 
but if so they have the wrong genus or the wrong number of punctures.
Note that the map $\tPar  \to \Par $, 
given by $(\nu; \Sigma) \mapsto \Sigma$,
is only outside of $\tPar'$ (respectively, $\Par'$) a covering.

\subsection{Vertical and horizontal boundary operator}

The boundary operator $\del = \del' + (-1)^q \del''$ 
on the complex $\bPar_\bubu$
is a sum of a vertical and a horizontal boundary operator, 
which are given by
\begin{align}
\del' &= \sum_{i=0}^q (-1)^i \del'_i \text{ with } 
\del'_i(\Sigma) = (\sigma_q, \ldots, \widehat{\sigma}_i, \ldots, \sigma_0) \\
\del'' &= \sum_{j=0}^p (-1)^j \del''_j \text{ with }
\del''_j(\Sigma) = (D_j(\sigma_q), \ldots, D_j(\sigma_0)).
\end{align}
Here $D_j \co \mf{S}_{p+1} \to \mf{S}_p$ deletes the letter $j$ 
from the cycle it occurs in and renormalizes the indices;
more precisely,
$D_j(\alpha) = s_j \circ \perm{\alpha(j),j} \circ \alpha \circ d_j$,
where $d_j \co [p-1] \to [p]$ is the simplicial degeneracy function
which avoids the value $j$,
and $s_i \co [p] \to [p-1]$ is the simplicial face function 
which repeats the value $i$.

The boundary operator $\tdel = \tdel' + (-1)^q \tdel''$ 
on the complex $\tbPar_\bubu$ has the following terms. 
The vertical face operators do not act on the new data $\nu$, 
thus
\begin{equation}
\tdel'_i(\nu;\Sigma) = (\nu;\sigma_q, \ldots, \widehat{\sigma}_i, \ldots, \sigma_0).
\end{equation}
For the horizontal face operators
we need a deletion operator $\Delta_j(\nu) = \nu \circ d_j$
defined by composing $\nu$ with the simplicial degeneracy function
$d_j \co [p-1] \to [p]$ which avoids the value $j$.
Then we define 
\begin{equation}
\tdel''_j(\nu;\Sigma) = (\Delta_j(\nu); D_j(\sigma_q), \ldots, D_j(\sigma_0)) .
\end{equation}

\subsection{Homogeneous and inhomogeneous notation}

It is convenient to have -- 
apart from the homogeneous notation above --
the inhomogeneous (or ``bar'') notation at hand: 
$\cT = [\tau_q | \ldots | \tau_1]$
where $\tau_k := \sigma_k \cdot \sigma_{k-1}^{-1}$  
is an element in $\fS_{p+1}$.
Vice versa, $\sigma_k := \tau_k \cdots \tau_1 \cdot \omega$.
The norm $\norm(\cT) = \ell(\tau_q) + \cdots + \ell(\tau_1) \leq h$,
and the number of cycles 
$\ncyc(\cT) = \ncyc(\sigma_q) = \ncyc(\tau_q \, \cdots \, \tau_1 \cdot \omega) \leq m+1$.
The conditions translate now into 
$\tau_k(0) = 0$ for $k = 1, \ldots, q$,
$\tau_k \not= \eins$ for $k = 1, \ldots, q$,
and that there is now common fixed point of all $\tau_k$, except $0$.

The vertical face operators are now
$\del'_i(\cT) = [\tau_q | \ldots |\tau_{i+1}\tau_i | \ldots]$,
for $i = 1, \ldots, q-1$, and are set equal to $0$ for $i=0,q$.
The horizontal face operators translate into 
corresponding operators $\del''_j(\cT)$.

In the non-permutable case, all conditions on the puncture enumeration $\nu$
must be expressed as before.


\section{Some properties of the chain complexes  $\bPar_\bubu$ and $\tbPar_\bubu$ }  
\label{sec5}

Let $h$, $m \geq 0$ be fixed, and set $n=1$.

\subsection{The complexes $\bPar_\bubu$}

The spaces $\Par$ (respectively, $\tPar$) are the geometric realizations 
of the bi-semi-simplicial complexes consisting of all the cells,
degenerate or not.

The degenerate cells form a subcomplex $\Par'$ resp $\tPar'$.
Thus they are finite cell complexes
whose cells are products $\Delta^p \times \Delta^q$ of simplices, 
where $p = 0,1, \ldots, 2h$ and $q = 0,1, \ldots, h$.
All cells of dimension less then $h+2$, in particular all $0$--cells,
are contained in $\Par'$ (respectively, in $\tPar'$).

We remark here that a subcomplex of $\Par$ is contractible.
Namely, if we change the definition of the norm by adding the term
$\norm(\sigma_0 \omega^{-1})$,
we have the same non-degenerate cells, but fewer degenerate cells.
The complement of these degenerate cells is still homeomorphic to the space $\Dip$.
After this modification we would have 
$H^{*-1}(\Par') \isomorphic H_{3h-*}(\Mod)$.
We will investigate this in a forthcoming article
exhibiting $\Mod$ as the interior of a compact manifold with boundary.

\subsection{The quotient complex} 

Since the subcomplex $\bPar'_\bubu$ 
is generated by basis elements,
the quotient complex 
$\bQ_\bubu := \bPar_\bubu/\bPar'_\bubu $ 
is generated by all $\Sigma$ not in the subcomplex; 
and the face operators $\del'_i$ (respectively, $\del''_j$) are set to be zero 
if they land in the subcomplex.
In particular $\del'_0 = \del'_q = 0$ and $\del''_0 = \del''_p = 0$,
since they always reduce the norm or the number of punctures
or violate any of the condition stated.

Similar statements hold for the chain complex 
$\tbQ_\bubu := \tbPar_\bubu/\tbPar'_\bubu$.

\subsection{The spectral sequence of the double complex}

The double complex $\bQ_\bubu$ gives rise to a 
spectral sequence starting with the $E^0$--term 
$E^0_{p,q} = \bQ_{p,q}$. 
It is concentrated in the first quadrant in the rectangle
$0 \leq p \leq 2h$ and $0 \leq q \leq h$; 
moreover, $E^0_{p,q} = 0$ for $p=0,1$ and for $p q \leq h$.

It turns out that (for fixed $p$) the vertical homology 
$E^1_{p,q} = H_q( \bQ_{p,\bullet}, \del')$
is concentrated in the maximal degree $q=h$.
Thus the $E_1$--term is a chain complex with differential
induced by $\del''$.
And the spectral sequence collapses with $E^2 = E^{\infty}$.

\subsection{The bar resolution}

For fixed $p$ the chain complex $\bPar_{p,\bullet}$ is similar to 
the homogeneous bar resolution of the group $\fS_p$.
First, since we require $\sigma_i(p) = 0$ for the permutations 
$\sigma_i \in \fS_{p+1}$, they are actually in one-to-one correspondence 
with $\fS_p$. Secondly, to normalize $\sigma_0 = \omega_p = \perm{0,1,\ldots,p}$
amounts to taking the coefficient module $\bZ$ with trivial $\fS_p$ action.
If one filters the bar resolution by the $\norm$, 
then $\bPar_{p,\bullet}(h,m,1)$ is a summand in the $h$th filtration quotient,
namely the summand determined by $\sigma_q$ having $m+1$ cycles.

Note that each cell $\Sigma$ intersects the subspace $\Par'$ of degenerate surfaces; 
more precisely, in the quotient chain complex 
$\bQ_{\bubu}$ we have
$\del'_0(\Sigma)  = \del'_q(\Sigma)   = 0$ and
$\del''_0(\Sigma) =  \del''_p(\Sigma) = 0$ for a cell $\Sigma$ of bi-degree $(p,q)$.

Similar statements hold for the chain complex $\tbQ_{p,\bullet}$. 

\section{The example $g=1$, $n=1$ and  $m=0$ algebraically}
\label{sec6}

In this section we discuss the example $g=1$, $m=0$ and $n=1$, the
moduli space of tori with one boundary curve and no punctures.  Thus $d
+ c = 6$.  The possible values for $p$ are $p=0,1,2,3,4$, and for $q$
they are $q=0,1,2$.  We write the cells in the homogeneous notation.

\subsection{The cells }

Recall the eight non-degenerate cells $\Sigma_1, \ldots, \Sigma_8$
as given above in \fullref{sec3.4}.
With bi-degrees
\begin{align*}
(p,q) &= (4,2) : \Sigma_1 \text{ and } \Sigma_2 \text{ in dimension } 6
\\
(p,q) &= (4,1) : \Sigma_3 \text{ in dimension } 5  \\
(p,q) &= (3,2) : \Sigma_4, \Sigma_5 \text{ and } \Sigma_6 \text{ in
dimension } 5  \\
(p,q) &= (3,1) : \Sigma_7 \text{ in dimension } 4 \\
(p,q) &= (2,2) : \Sigma_8 \text{ in dimension } 4
\end{align*}

\subsection{The $E^0$--term}

The $E^0$--term is shown in the following table:
\begin{equation}
\begin{array}{|c|c|c|c|c|c|} \hline
q=2 &     &     & \Sigma_8 & \Sigma_4, \Sigma_5, \Sigma_6 & \Sigma_1, \Sigma_2   \\ \hline
q=1 &     &     &          & \Sigma_7                     & \Sigma_3             \\ \hline
q=0 &     &     &          &                              &                      \\ \hline
    & p=0 & p=1 & p=2      & p=3                          & p=4
\\ \hline
\end{array}
\end{equation}
The vertical boundary operator $\partial'$ is
$$\partial'(\Sigma_1) = \partial'(\Sigma_2) = -\Sigma_3,\qquad
\partial'(\Sigma_4) = \partial'(\Sigma_5) = \partial'(\Sigma_6) =
-\Sigma_7.$$

\subsection{The $E^1$--term}

We set $A := \Sigma_1 - \Sigma_2$, 
$B := \Sigma_4 - \Sigma_5$, and $C := \Sigma_4 - \Sigma_6$.
Note that in this example the fundamental class
(see \fullref{sec7}) is $\mu = -A$.

Then the $E^1$--term is as follows:
\begin{equation}
\begin{array}{|c|c|c|c|c|c|}                 \hline
q=2 &     &     & \Sigma_8 & B, C   & A   \\ \hline 
q=1 &     &     &          & 0      & 0   \\ \hline
q=0 &     &     &          &        &     \\ \hline
    & p=0 & p=1 & p=2      & p=3    & p=4 \\ \hline
\end{array}
\end{equation}
The horizontal boundary operator $\partial''$ is
$$\partial''(A) = 0,\qquad
\partial''(B) = -2\Sigma_8,\qquad
\partial''(C) = -\Sigma_8.$$

\subsection{The $E^2$--term and final result}

If we set $D := B - 2C$ we have a cycle in degree $(3,2)$,
and therefore an $E^2$--term as follows:
\begin{equation}
\begin{array}{|c|c|c|c|c|c|}            \hline
q=2 &     &     &  0  &  D   &  A    \\ \hline
q=1 &     &     &     &  0   &  0    \\ \hline
q=0 &     &     &     &      &       \\ \hline
    & p=0 & p=1 & p=2 & p=3  & p=4   \\ \hline
\end{array}
\end{equation}
The moduli space of tori with one boundary curve is
orientable and homotopy-equivalent to the complement of the trefoil
knot in $\bR^3$, or in other words 
the mapping class group $\Gamma$ is isomorphic to the 
third braid group; its cohomology is infinite cyclic in degrees 
$0$ and $1$.


\section{The fundamental class}
\label{sec7}

Here we set $n=1$. Then $h = 2g+m$.
The (relative) cycle $\mu$ representing the fundamental class $[\mu] \in H_{3h}(\bPar,\bPar';\bZ)$ 
is given as follows. 
First, its degree is $(p,q) = (2h,h)$.
Thus we need $2h$ permutations $\tau_h, \ldots, \tau_1$ for the inhomogeneous notation.
We set $\tau_1 = \perm{3,1}$ and $\tau_2 = \perm{4,2}$,
an interlocking pair of disjoint transpositions.
We continue in this way until we have $g$ blocks of two pairs of interlocking 
transpositions.
Then we continue with $\tau_{2g+1} = \perm{4g+2,4g+1}$,
till $\tau_{2g+m} = \perm{4g+2m,4g+2m-1}$;
these are $m$ neighbour transpositions.

Altogether we have $2h$ disjoint (and therefore commuting) transpositions.
The $g$ {\it symplectic pairs} of interlocking transpositions 
create each a handle of the surface; 
and the $m$ non-interlocking transpositions 
create each a puncture.
The cycle number of $\tau_h \cdots \tau_1 \cdot \omega_{2h}$ is $m+1$.

Let $\Lambda(h,m)$ denote the subset of all permutations $\alpha \in \fS_{2h}$
such that the horizontally conjugated cell
\begin{equation}
\kappa_\alpha (T) = [\alpha \tau_h \alpha^{-1} | \ldots | \alpha \tau_1 \alpha^{-1} ]
\end{equation}
has the same cycle number $\ncyc(\alpha . T) = m+1$.
Furthermore, for $\pi \in \fS_h$ let   
\begin{equation}
\pi . T = [\tau_{\pi(h)} | \ldots | \tau_{\pi(1)} ]
\end{equation}
denote the vertically permuted cell.
These operations were studied by Schardt \cite{Schardt}.
The chain  
\begin{equation}
\mu = 
\sum_{\pi \in \fS_{h}} \,\,\, \sum_{\alpha \in \Lambda(h,m)} 
\sign(\pi) \sign(\alpha) \,\, \pi \,.\, \kappa_\alpha(T)
\end{equation}
is a cycle and represents the fundamental class.
For example, when $g=2, m=0$, this cycle has $504$ terms.

\section{The orientation system}
\label{sec8}

The local coefficient system is the orientation system of the relative manifold
$\Par, \Par'$. For $m \geq 2$ it is not trivial, see M\"uller \cite{Mueller}.
In order to determine the change of sign along a path,
we can restrict to edge paths from barycenters to barycenters in a subdivision.

We will associate to each cell $\cT$ a top-dimensional cell $\topcell(\cT)$.
Denote then by $\epsilon(\cT)$ 
the sign with which $\topcell(\cT)$ occurs in the fundamental class $\mu$.
The sign change for the path from the barycenter of $\cT$ 
to the barycenter of $\cT'$, where $\cT' \subset \cT$, 
is then the product $\epsilon(\cT) \epsilon(\cT')$ of these two signs.

To define the cell $\topcell(\cT)$ we associate to $\cT = [\tau_q | \ldots | \tau_1]$
we start by writing each $\tau_i$ in a certain normal form $\mathrm{NF}(\tau_i)$.
First, if $\alpha = \perm{i_l, \ldots, i_0}$ is a permutation 
consisting of just one cycle,
we write it as 
$\alpha = \mathrm{NF}(\alpha) = \perm{i_l,i_{l-1}} \cdot \cdots \cdot \perm{i_1,i_0}$.
If $\alpha$ consists of several (non-trivial) cycles, 
we order them by their minima, the absolute minimum being first 
(or in our notation rightmost).
This gives a normal form
\begin{equation}
\alpha = \mathrm{NF}(\alpha) = 
\cdots 
       \cdot \perm{i'_{l'},i'_{l'-1}} \cdot \cdots \cdot \perm{i'_1,i'_0} 
       \cdot \perm{i_l,i_{l-1}}       \cdot \cdots \cdot \perm{i_1,i_0}.
\end{equation}
Assume we have the normal forms $\mathrm{NF}(\tau_q), \ldots, \mathrm{NF}(\tau_1)$. 
The next step is to ``separate'' the factors in each normal form $NF(\tau_i)$
by replacing all product signs by a bar symbol.
We obtain a new cell:
\begin{equation}
\cT' = [\theta_h \,|\, \theta_{h-1} \,|\, \ldots \,|\, \theta_1 ].
\end{equation}
Note that it has exactly $h$ transpositions $\theta_h, \ldots, \theta_1$. 
Thus $\cT$ has bi-degree $(p,h)$ 
(and is not yet the top-dimensional cell we seek). 

The next step is to ``spread out'' these transpositions, that is
to replicate indices if they occur in several of the transpositions.
Denote for $j=0, \ldots, p$ by $S_j \co \fS_p \to \fS_{p+1}$ the function
which sends the permutation $\alpha$ to the new permutation $S_j(\alpha)$
by increasing the indices $k \geq j$ in a cycle notation to $k+1$
and leaving the others. Note that $j$ becomes a fixed point.  In other
words, $S_j(\alpha) = d_j \circ \perm{\alpha(j),j} \circ \alpha \circ
s_j$ expressed with the simplicial face and degeneracy maps.  (See the
definition of $D_j$ in \fullref{sec4.2}; note that neither the $D_j$
nor the $S_j$ are group homomorphisms; they turn the family of symmetric
groups into a crossed simplicial group.)  Since $S_j(\perm{k,l}) =
\perm{d_j(k), d_j(l)}$, a transposition is mapped to a transposition.

Assume $k$ is the lowest index in $\cT'$ occurring more than ones in one
of transpositions $\theta_i$, say in $\theta_{i_1}$ the first time and
in $\theta_{i_2}$ the second time, for $i_2 > i_1$.  Then we replace $\cT'$ by
\begin{equation}
\cT'' =
[  S_k(\theta_h)        |  S_k(\theta_{h-1})  | {\ldots}   |
   S_k(\theta_{i_2})     |                          {\ldots}   |
  S_k(\theta_{i_{1} +1}) | 
   S_{k+1}(\theta_{i_1})  |                         {\ldots}   | S_{k+1}(\theta_1)            ],
\end{equation}
a cell of bi-degree $(p+1,h)$.  We repeat this last step until no index
occurs more than ones in these transpositions, each time raising the
first degree by one.  The final result is a sequence of $h$ disjoint
transpositions, thus a cell $\topcell(\cT)$ of degree $(2h,h)$.

\section{Comments on the computation}
\label{sec9}

The C++ program with documentation is contained in the first author's
thesis \cite{Abhau}.  It is about 3000 lines long and consists of 187
subroutines.  Here we give some comments.

\subsection{Listing the cells}

The computer program uses the C++ Standard Template Library extensively. 
A permutation in $\fS_{p+1}$ is implemented as a vector of integers (list of values)
and also as a single integer using an easy computable bijection 
$\mathfrak{S}_{p+1} \rightarrow \{1,\ldots, (p+1)!\}$.
Thus manipulations can be performed easily and an efficient 
storage handling is possible up to $h = 6$. 

We obtain the non-degenerate cells or generators of $\bQ_{\bubu}$ 
by first generating the set $\Lambda$ of all cells in bi-degree $(2h,h)$,
namely all $h$--tuples of disjoint transpositions in $\fS_{2h}$
-- there are $L_h = (2h-1)!! = \frac{(2h)!}{2^h h!}$ such tuples --
then we sort them according to their cycle number $m$ into sets $\Lambda(h,m)$.
Note that for $h$ odd (respectively, even) the cycle number $m$ 
is again odd (respectively, even).

Applying the vertical boundary operator $\del'$ successively 
gives all generators of the chain complex $\bQ_{2h,\bullet}$. 
Next we generate the cells of bi-degree $(2h-1,h)$ 
as images of the cells of bi-degree $(2h,h)$ 
under the horizontal boundary operator $\del''$. 
Inductively, we obtain all cells of $\bQ_{\bubu}$. 

The same procedure gives the cells of the complex 
$\tbQ_{\bubu}$, 
where we store all possible puncture enumerations  
$\nu \co [p] \rightarrow [m]$ additionally.

\subsection{The matrices of the boundary operator}

While generating the cells, 
the matrices $D'$ resp $D''$ of the two boundary operators 
$\del'$ (respectively, $\del''$) are produced simultaneously. 
We store the matrices in the following sparse format: 
for the matrix entries, a compressed row format is used;
moreover, a compressed column matrix of boolean variables gives the positions 
of the nonzero entries. 
This seems to be a reasonable trade-off between memory space and running time, 
because both row and column operations are performed frequently on the matrix.

For the homology of the complex $\bQ_{\bubu}$ 
(respectively, $\tbQ_{\bubu}$)
we first compute the Smith normal forms of the matrices 
of the vertical boundary operators. 
As noted earlier, for each $p$, the homology of the vertical chain complexes 
$\bQ_{p,\bullet}$ is just the kernel of 
$\del' \co \bQ_{p,h} \to \bQ_{p,h-1}$.
This gives the $E^1$--term of the spectral sequence.

In order to compute its $E_2$--term and thus the desired homology 
of the total complex,
we add a lattice reduction routine, the LLL-algorithm 
(named after Lenstra-Lenstra-Lov{\'a}sz) 
while computing the Smith normal form (or elementary divisors)
of the horizontal boundary operator $\del''$.
This is important to avoid a coefficient explosion in the case $h=5$. 
All the computations were also done modulo 2.

\subsection{Size of the matrices}

Both running time and memory space grow exponentially with $h$. 
More precisely, let $N_{p,q} = N_{p,q}(h,m)$ be the number of cells of bi-degree $(p,q)$.
\fullref{tbl1} shows these numbers for some cases.

\subsubsection*{Ranks of the $E^0$--term for $h=5$ and $q=5,4$}

\begin{table}[ht!]\tabcolsep 4.5pt\centering
\begin{small}
\begin{tabular}{|c|c|c|c|c|c|c|c|c|c|c|c|}\hline
$p {=}$            & 0   & 1   & 2   & 3     & 4      & 5       & 6        & 7        & 8        & 9        & 10      \\ \hline
\begin{tabular}{|c|} \hline $h {=} 5$ \\ $m {=} 1$ \\ \hline \end{tabular}
   &     &     &     &       &        &         &          &          &          &          &         \\ 
\hfill     $q{=}5$ & $0$ & $0$ & $1$ & $240$ & $6170$ & $51115$ & $195264$ & $394240$ & $435680$ & $249480$ & $57960$ \\       
\hfill     $q{=}4$ & $0$ & $0$ & $0$ & $216$ & $7840$ & $76140$ & $320880$ & $694148$ & $808192$ & $482328$ & $115920$ \\ \hline
%
\begin{tabular}{|c|} \hline $h {=} 5$ \\ $m {=} 3$ \\ \hline \end{tabular}
   &     &     &     &       &        &         &          &          &          &          &         \\ 
\hfill     $q{=}5$ & $0$ & $0$ & $0$ & $0$   & $640$  & $12425$ & $74610$  & $202825$ & $278600$ & $189000$ & $50400$ \\       
\hfill     $q{=}4$ & $0$ & $0$ & $0$ & $0$   & $800$  & $18500$ & $122700$ & $357280$ & $516880$ & $365400$ & $100800$ \\ \hline
%
\begin{tabular}{|c|}
\hline $h {=} 5$ \\ $m {=} 5$ \\ \hline \end{tabular}
   &     &     &     &       &        &         &          &          &          &          &         \\ 
\hfill     $q{=}5$ & $0$ & $0$ & $0$ & $0$   & $0$    & $0$     & $1296$   & $7735$   & $16520$  & $15120$  & $5040$  \\       
\hfill     $q{=}4$ & $0$ & $0$ & $0$ & $0$   & $0$    & $0$     & $2160$   & $13692$  & $30688$  & $29232$ & $10080$  \\ \hline
\end{tabular}
\end{small}
\caption{Ranks of the $E^0$--term for $h=5$ and $q=5,4$}
\label{tbl1}
\end{table}

Observe that both boundary operators leave the puncture number $m$ invariant, 
so we can compute the homology for each
$m \equiv h \mod 2$ separately.  

To generate the matrix $D'_{p,q}$ or $D''_{p,q}$ of a boundary operator
$\del' \co \bQ_{p,q} \to \bQ_{p,q-1}$ (respectively,
$\del'' \co \bQ_{p,q} \to \bQ_{p-1,q}$)
we need to take $q-1$ (respectively, $p-1$) boundary faces from $N_{p,q}$ cells 
and search these face cells in a (sorted) list of $N_{p,q-1}$ resp $N_{p-1,q}$ cells. 
This leads to a running time of $O( h N^2 \log N )$, 
where $N = \max \{N_{p,q} \}$. 
We only need to compute the matrices $D'_{2,h}, \ldots, D'_{2h,h}$ 
and $D''_{3,h}, \ldots, D''_{2h,h}$. 
See \fullref{tbl1} for the size of the matrices $D'_{p,5}$ for $h=5$.

The most time-consuming procedure 
when determining the Smith normal form is the LLL-algorithm. 
This algorithm, applied to $r$ columns of an integer matrix 
takes $O(r^6 \log^3 B)$ steps, 
where $B$ is the maximal norm of a column. 
But note that we apply the LLL-algorithm only 
when the coefficients exceed a prescribed barrier, 
and also only to portions of the columns, 
so $r$ is much smaller than $N$.

\subsubsection*{Ranks of the $E^1$--term for $h=5$}

The most memory-consuming part in the computation procedure is the storage
of the matrices $D''_{p,h}$ for the horizontal boundary operator $\del''$.
in terms of the transformed basis.  Numerical values for the number
of generators of the $E_1$--term can again be taken from \fullref{tbl2}
giving estimates for the matrix sizes $D''_{p,5}$.

\begin{table}[ht!]\centering
\begin{small}
\begin{tabular}{|r|c|c|c|c|c|c|c|c|c|c|c|} \hline
$p =$            & 0   & 1   & 2   & 3     & 4      & 5       & 6        & 7        & 8        & 9        & 10      \\ \hline
\begin{tabular}{|c|} \hline    $h = 5$ \\ $m = 1$ \\ \hline \end{tabular}
   &     &     &     &       &        &         &          &          &          &          &         \\ 
\hfill     $q=5$ & $0$ & $0$ & $1$ & $60$  & $650$  & $2860$  & $6588$   & $8708$   & $6678$   & $2772$   & $483$   \\ \hline   
\begin{tabular}{|c|} \hline    $h = 5$ \\ $m = 3$ \\ \hline \end{tabular}
   &     &     &     &       &        &         &          &          &          &          &         \\ 
\hfill     $q=5$ & $0$ & $0$ & $0$ & $0$   & $70$   & $700$   & $2520$   & $4480$   & $4270$   & $2100$   & $420$   \\ \hline   
\begin{tabular}{|c|} \hline    $h = 5$ \\ $m = 5$ \\ \hline \end{tabular}
   &     &     &     &       &        &         &          &          &          &          &         \\ 
\hfill     $q=5$ & $0$ & $0$ & $0$ & $0$   & $0$    & $0$     & $1$      & $14$     & $56$     & $84$     & $42$    \\ \hline  
\end{tabular}
\smallskip
\end{small}
\caption{Ranks of the $E^1$--term for $h=5$}
\label{tbl2}
\end{table}

We were able to compute the homology of the complexes $\bQ_{\bubu}$
up to $h=5$ on a machine with 4 GB of RAM and a 3 GHz processor taking
approximately 4 days for the case $h=5$, $m=1$.  The computational
problems for $h > 5$ are both running time and memory space.

\bibliographystyle{gtart}
\bibliography{link}

\end{document}